\documentclass[a4paper]{article}

\usepackage{amssymb, amsmath}
\usepackage{amsthm}

\newtheorem{lemma}{\bf Lemma}[section]%
\newtheorem{theorem}[lemma]{\bf Theorem}%
\newtheorem{proposition}[lemma]{\bf Proposition}% 
\newtheorem{corollary}[lemma]{\bf Corollary}%   
\newtheorem{example}[lemma]{\bf Example}%    
\newtheorem{definition}[lemma]{\bf Definition}{}%  
\newcommand{\be}{\begin{equation}}
\newcommand{\ee}{\end{equation}}
\newcommand{\ol}{\overline}
\newcommand{\wt}{\widetilde}
\newcommand{\wh}{\widehat}

 \title{Graph polynomials and Tutte-Grothendieck invariants: \\ an application of elementary finite Fourier analysis}
 \author{Andrew Goodall}
\begin{document}
\maketitle
\begin{abstract}
This paper is based on a series of talks given at the Patejdlovka
Enumeration Workshop held in the Czech Republic in November 2007. The
topics covered are as follows. The graph polynomial, Tutte-Grothendieck invariants, an overview of relevant elementary finite Fourier analysis, the Tutte polynomial
of a graph
as a Hamming weight enumerator of its set of tensions (or flows), description
of a family of
polynomials containing the graph polynomial which yield
Tutte-Grothendieck invariants in a similar way.
\end{abstract}

\section{Introduction}

The graph polynomial is a generalization of the Vandermonde
determinant (which may be viewed as the graph polynomial of a complete graph) that
was considered by Sylvester and Petersen in the nineteenth century in
connection with binary quantic forms. Alon,
Tarsi and Matiyasevich in more recent years have found that it contains
a good deal of information about the vertex colourings of a graph. For
example, the number of proper $3$-colourings of a graph
is a simply defined function of the coefficients of its graph
polynomial. In this
article we consider a family of polynomials containing the graph
polynomial and ask whether other Tutte-Grothendieck invariants
can be obtained in a similar way. Our results are obtained by expressing the relevant
parameter as
the partition function of a vertex colouring model (such as the Potts
model) or, in different language, the graph parameter obtained from homomorphisms to a weighted graph.  

In Section \ref{sec: Petersen graph poly} we define the graph
polynomial and explore its relation to proper
vertex colourings. In Section \ref{sec: TG} Tutte-Grothendieck
invariants are defined and their pervasiveness noted. In Section
\ref{sec: Fourier} a potted account is given of
Fourier analysis on finite Abelian groups which will be used to obtain our results. In Section
\ref{sec: we} tensions
and flows of a graph are defined and the view of the Tutte polynomial as a Hamming
weight enumerator propounded. In the final Section \ref{sec: graph polynomials} we
characterize those polynomials which share with the graph
polynomial the property of yielding a Tutte-Grothendieck invariant from the $\ell_2$-norm of their coefficients. More generally, the graph
polynomial is seen to belong to a family of polynomials
 whose $\ell_2$-norm is equal to an
 evaluation of the complete weight enumerator of the
set of tensions (or flows) of the graph. 

An expanded version of Section \ref{sec: Fourier} can be found in
\cite{AJG07}, and an even more fulsome presentation is given in \cite{CCC07}. 
The book \cite{Terras99} is recommended for an introduction to finite
Fourier analysis and its wide range of applications.

\section{The graph polynomial}\label{sec: Petersen graph poly}
Let $G=(V,E)$ be a graph 
with some fixed, arbitrary
orientation of its edges, and denote its directed edge set by
$\overrightarrow{E}$.

Let $Q$ be a finite set of size $q$. A
{\em proper vertex $q$-colouring} using colour set $Q$ is an assignment of colours $(c_v:v\in
V)\in Q^V$ such that $c_u\neq c_v$ whenever $\{u,v\}\in E$. 
The number of proper vertex $q$-colourings of $G$ is denoted by
$P(G;q)$ (an
evaluation of the chromatic polynomial of $G$ at $q$). 

Let $\mathbf{x}=(x_v: v\in V)$ be a tuple of commuting indeterminates
indexed by $V$ and
define the {\em graph polynomial}\,\footnote{The graph polynomial has
  not yet acquired the qualification of a proper name. 
The `Sylvester-Petersen polynomial' might be a candidate
\cite{Sylvester1878, P1891}. Matiyasevich analyses the graph polynomial of the line
graph of a cubic plane graph in order to obtain reformulations of the
Four Colour Theorem \cite{M04}. Alon and Tarsi \cite{AT92, AT97,
  Tarsi} interpret its
coefficients in terms of orientations; their interpretations in
terms of proper vertex colourings will be described in this
section. Ellingham  and Goddyn \cite{EG96} call the graph polynomial
the {\em graph  monomial} averring that the latter has a less
anonymous character than the former.} $F(G)$ in $\mathbb{C}[\mathbf{x}]$ by
$$F(G;\mathbf{x})=\prod_{(u,v)\in\overrightarrow{E}}(x_u-x_v).$$
Given an assignment of values $\mathbf{c}=(c_v:v\in V)\in\mathbb{C}^V$
to the indeterminates $\mathbf{x}=(x_v:v\in V)$, the graph polynomial takes a non-zero value
if and only if
$\mathbf{c}$ corresponds to a proper vertex colouring with colour set
$Q=\{c_v:v\in V\}$. Set $\zeta=e^{2\pi i/q}$. By restricting $c_v$ to one of the $q$ points
$1,\zeta,\ldots,\zeta^{q-1}$ on the unit circle  a criterion emerges
for the existence of a proper vertex $q$-colouring of $G$ in terms of
the polynomial $F(G;\mathbf{x})$.
 
The algebraic
variety of points $\{(c_v:v\in V):c_v\in
\{1,\zeta,\ldots,\zeta^{q-1}\}\}$ corresponds to the ideal
$(x_v^q-1:v\in V)$ of the ring $\mathbb{C}[\mathbf{x}]$.
Denote the graph polynomial modulo the ideal generated by the
polynomials $x_v^q-1$ by
$$F^{(q)}(G;\mathbf{x})=F(G;\mathbf{x})\bmod (x_v^q-1:v\in V).$$

By Lagrange interpolation,
\begin{align*} F^{(q)}(G;\mathbf{x}) & =\sum_{(a_v:v\in V)}\prod_{v\in V}\prod_{a\neq
  a_v}\frac{x_v-\zeta^{a}}{\zeta^{a_v}-\zeta^{a}}\;F(G;(\zeta^{a_v}:v\in V)\,)\\
 & = q^{-|V|}\sum_{(a_v:v\in V)}\prod_{v\in
  V}\frac{x_v^q-1}{\zeta^{-a_v}x_v-1}F(G;(\zeta^{a_v}:v\in V)\,),\end{align*}
 where the summations are over $(a_v:v\in
V)\in\{0,1,\ldots, q-1\}^V$ and the last line follows since $\prod_{a\neq a_v}(\zeta^{a_v}-\zeta^{a})=\zeta^{(q-1)a_v}\prod_{b\neq 0}(1-\zeta^b)=\zeta^{-a_v}q.$
The relationship between coefficients of the polynomial $F^{(q)}(G;\mathbf{x})$ and
its evaluations at points $(\zeta^{a_v}:v\in V)$ is exhibited here as
a basis change between the basis of monomials $\prod_{v\in
  V}x_v^{a_v}$ and the basis of polynomials $\prod_{v\in
  V}\frac{x_v^q-1}{\zeta^{-a_v}x_v-1}$. The connection is the Fourier transform. This article is
an elaboration of this remark.

Alon and Tarsi \cite{AT97} use the ``Combinatorial Nullstellensatz'' \cite{A99}
to prove that
$F^{(q)}(G;\mathbf{x})\neq 0$ if and only if
$P(G;q)>0$, and also show that more can be said.

For a polynomial
$F(\mathbf{x})\in\mathbb{C}[\mathbf{x}]/(x_v^q-1:v\in V)$, let
$[\mathbf{x}^\mathbf{a}]F(\mathbf{x})$ denote the coefficient of the
monomial $\mathbf{x}^{\mathbf{a}}=\prod_{v\in V}x_v^{a_v}$ in its
expansion to the monomial basis. In particular,
$[\mathbf{x}^\mathbf{0}]F(\mathbf{x})$ is the constant term of $F(\mathbf{x})$. The (squared) $\ell_2$-norm
$\|F(\mathbf{x})\|_2^2$ is defined by
$$\|F(\mathbf{x})\|_2^2=\sum_{\mathbf{a}\in\{0,1,\ldots,q-1\}^V}\big|[\mathbf{x}^\mathbf{a}]F(\mathbf{x})\big|^2.$$
% to be the sum of the absolute squares of its  coefficients. 
That this is a norm includes the fact that
  $F^{(q)}(G;\mathbf{x})\neq 0$ if and only if
  $\|F^{(q)}(G;\mathbf{x})\|_2^2\neq 0$.

\begin{theorem}\label{thm: Alon Tarsi} {\rm  \cite{AT97}} For each $q\in\mathbb{N}$,
 $$\|F^{(q)}(G;\mathbf{x})\|_2^2=q^{-|V|}4^{|E|}\sum_{\mathbf{c}\in\{0,1,\ldots,
    q-1\}^V}\prod_{uv\in E}\sin^2\frac{\pi(c_v-c_u)}{q},$$
 the sum being over all vertex colourings of $G$ with colours
 $\{0,1,\ldots, q-1\}$. In particular,  for $q=3$ this is $3^{|E|-|V|}P(G;3)$.
\end{theorem}

For the next theorem we require a further definition. A {\em $(q,1)$-flow} of
$G$ is a partial orientation of $G$ with the property that at each
vertex the number of edges directed out of $v$ is congruent to the
number of edges directed into $v$ modulo $q$. (A partial orientation
is obtained when some edges of $G$ are assigned an orientation while
the other edges remain undirected.) By referring to the fixed
orientation $\overrightarrow{E}$ of $G$, it is possible to
use the equivalent definition as an assignment of values $(b_e:e\in
E)$ to the edges of $G$ with the properties that  $b_e\in\{0,1,-1\}$ 
and the net flow (incoming minus outgoing values) at each vertex is equal to
zero modulo $q$.

\begin{theorem}\label{thm: Tarsi} {\rm \cite{Tarsi}} 
For each $q\in\mathbb{N}$,
 $$\|F^{(q)}(G;\mathbf{x})\|_2^2=(-1)^{|E|}\sum_{(q,1)\mbox{\rm \small -flows}\;
  \mathbf{b}}(-2)^{|E|-|\mathbf{b}|},$$
where
$|\mathbf{b}|=\#\{e\in E:b_e\neq 0\}$. 
\end{theorem}

One aim of this article is to reveal the underlying relationship
between  Theorems \ref{thm: Alon Tarsi} and \ref{thm: Tarsi} in a more
general context. The other is to characterize those polynomials of the form
$$\prod_{(u,v)\in\overrightarrow{E}}\;\sum_{a,b\in\{0,1,\ldots,
  q-1\}}\;\;f(a,b)x_u^ax_v^b\hspace{1cm}\mbox{\rm
  mod}\, (x_v^q-1:v\in V)$$
whose $\ell_2$-norm is a Tutte-Grothendieck invariant (such as
$P(G;q)$). The graph polynomial is the case $f(1,0)=1, f(0,1)=-1$ and
$f(a,b)=0$ otherwise, and Theorem \ref{thm: Alon Tarsi} says that for
$q=3$ its $\ell_2$-norm is the Tutte-Grothendieck invariant $3^{|E|-|V|}P(G;3)$.

%In the next section we define Tutte-Grothendieck
%invariants in order to be able to begin searching for them amongst
%polynomials such as the graph polynomial. 

\section{Tutte-Grothendieck invariants}\label{sec: TG}

Let $G=(V,E)$ be a graph, loops and parallel edges permitted,  with
$k(G)$ components, rank $r(G)=|V|-k(G)$ and nullity $n(G)=|E|-r(G)$.
 
{\em Deleting} an edge $e\in E$ gives a graph $G\setminus e$ with one fewer
edge than $G$. {\em Contracting} $e$ gives a graph $G/e$ with one fewer vertex
and one fewer edge than $G$. Many graph parameters may be
recursively defined via contraction-deletion recurrences.

\begin{definition}
A function $F$ from (isomorphism classes of) graphs to $\mathbb{C}[\alpha,\beta,\gamma,x,y]$ is a {\em Tutte-Grothendieck
  invariant} if it satisfies, for each graph $G=(V,E)$ and any edge
$e\in E$,
\be\label{TG}F(G)=\begin{cases} \gamma^{|V|} & E=\emptyset,\\
 xF(G/e) & e \;\; \mbox{\rm a bridge,}\\
 yF(G\setminus e) & e \;\;\mbox{\rm a loop,}\\
 \alpha F(G/e)+\beta F(G\setminus e) & e \;\;\mbox{\rm not a bridge or loop.}\end{cases}
\ee 
\end{definition}

See for example the accounts in
\cite{DW93, BB98, GR01} for an appreciation of the ubiquity of Tutte-Grothendieck
invariants.
For $A\subseteq E$, the subgraph $(V,A)$ is obtained from $G$ by
deleting edges not in $A$. Given $G=(V,E)$, the rank of the graph
$(V,A)$ is denoted by $r(A)$.
A Tutte-Grothendieck invariant is an evaluation of the {\em Tutte
  polynomial}, defined by
\be\label{eqn: Tutte}T(G;x,y)=\sum_{A\subseteq E}(x-1)^{r(E)-r(A)}(y-1)^{|A|-r(A)}.\ee
The Tutte polynomial is a rescaling of the Whitney rank polynomial of $G$
(for which see for example \cite{GR01}), a generating function for $(|A|,r(A))$ over all
subgraphs $(V,A)$ of $G$. The coefficients of the Tutte polynomial are
non-negative integers (see for example \cite{B94,BB98}), a fact while not evident from its definition in
equation \eqref{eqn: Tutte} is more readily seen in its alternative
formulation as a Tutte-Grothendieck invariant with $\alpha=\beta=\gamma=1$.

\begin{theorem}\label{thm: TG}
 If $F$ is a Tutte-Grothendieck invariant satisfying the equations
 \eqref{TG} then 
$$F(G)=\gamma^{k(G)}\alpha^{r(G)}\beta^{n(G)}T(G;\frac{x}{\alpha},\frac{y}{\beta}).$$
\end{theorem}

See \cite{BB98} for how to interpret this evaluation for the case $\alpha=0$ or $\beta=0$.

\begin{example}\label{ex: monochrome} The {\em monochromial} $P(G)=P(G;q,y)$ (bad colouring
    polynomial, coboundary polynomial, partition function of the $q$-state Potts model) is defined by 
\be\label{eqn: monochrome poly} P(G;q,y)=\sum_{\mathbf{c}\in
    Q^V}y^{\#\{(u,v)\in\overrightarrow{E}:c_u=c_v\}},\ee where $Q$ is
a set of $q$ colours (states) and $\mathbf{c}=(c_v:v\in V)$ is a
vertex colouring of $G$ using colours from $Q$.
It is easily verified that the function $P$ satisfies
$$P(G)=\begin{cases} q^{|V|} & E=\emptyset,\\
(y+q-1)P(G/e) & \mbox{\rm $e$ a bridge,}\\
yP(G\setminus e) & \mbox{\rm $e$ a loop,}\\
(y-1)P(G/e)+P(G\setminus e) & \mbox{\rm $e$ not a bridge or
  loop.}\end{cases}$$
By Theorem \ref{thm: TG}, 
\be\label{eqn: mono poly}P(G;q,y)=q^{k(G)}(y-1)^{r(G)}T(G;\frac{y\!-\!1\!+\!q}{y\!-\!1},y).\ee
In particular, the {\em chromatic polynomial} $P(G;q)$,
counting the number of proper vertex $q$-colourings of $G$, is given by
$$P(G;q)=q^{k(G)}(-1)^{r(G)}T(G;1-q,0).$$
\end{example}
 
Let $Q$ be a set of size $q$ (later in this article to be an additive Abelian group of order $q$) and
 $\mathbf{w}=(w_{a,b})$ a tuple of complex numbers indexed by
$(a,b)\in Q\times Q$. Assume that the edges $\{u,v\}$ of $G=(V,E)$ have been
given an arbitrary, fixed orientation $(u,v)$. Denote by
$\overrightarrow{E}$ the resulting set of directed edges. Consider the partition function for a vertex $Q$-colouring model that assigns a weight $w_{a,b}$ to a directed edge $(u,v)$ coloured $(a,b)$:
\begin{equation}\label{eqn: S}
F(G;\mathbf{w})=\sum_{\mathbf{c}\in Q^V}\prod_{(u,v)\in\overrightarrow{E}}w_{c_u,c_v}=\sum_{\mathbf{c}\in Q^V}\prod_{(a,b)\in Q\times Q}w_{a,b}^{\;\;\;\;\#\{(u,v)\in\overrightarrow{E}\,:\,(c_u,c_v)=(a,b)\}}.
\end{equation}
This partition function may be interpreted as the weight
of a graph homomorphism $G\rightarrow H$, where $H$ is the directed graph on
vertex set $Q$ and edge set $\{(a,b):w_{a,b}\neq 0\}$, with edge
weights $w_{a,b}$, i.e., the weighted graph $H$ has adjacency matrix $(w_{a,b})_{a,b\in Q}$.  %[and vertex weights $y_a$ ....]. 
(It is possible to also have vertex weights for $H$ in addition to its edge
weights, but this will not be considered here. See for example
\cite{FLS04, S05} and \cite{HN04, GNR07} for more on vertex
colouring models and on graph homomorphisms.)

\begin{theorem}\label{TG Hamming}
The graphical invariant $F(G;\mathbf{w})$ defined by equation \eqref{eqn: S} is a Tutte-Grothendieck invariant if and only if there are constants $y,w$ such that 
$$w_{a,b}=\begin{cases} w & a\neq b,\\
 y & a=b.\end{cases}
$$ 
In this case $F(G;\mathbf{w})=F(G;w,y)=q^{k(G)}w^{n(G)}(y-w)^{r(G)}T(G;\frac{y-(q\!-\!1)w}{y-w},\frac{y}{w})$. (If $w=0$ then $F(G;0,y)=y^{|E|}$ and if $w=y$ then $F(G;y,y)=q^{|V|}y^{|E|}$.)
\end{theorem}

A sketch only of a proof of Theorem \ref{TG Hamming} is
given.\footnote{I am grateful to Delia Garijo for alerting me to the
  fact that I was assuming the truth of something that required
  proof, and also for her description of how she has been tackling a related,
  stronger result.} 
%The
%reader is invited to feel some relief that it is but
%a sketch, for its details are somewhat wearing. 
(A fuller version will
appear in a forthcoming paper.) The
following lemma is the main tool.

\begin{lemma}\label{lemma: mth powers}
If $u_1,\ldots, u_r,v\in\mathbb{C}$ satisfy
\begin{equation}\label{eqn: all m} u_1^m+u_2^m+\cdots +u_r^m=rv^m,\ee
for all integers $m\geq 0$, then 
$$u_1=u_2=\cdots =u_r=v.$$
\end{lemma}
\begin{proof}
We may assume $v\neq 0$. Rewriting the relation \eqref{eqn: all m} in terms of ordinary
generating functions, 
$$(1-u_1z)^{-1}+(1-u_2z)^{-1}+\cdots + (1-u_rz)^{-1}=r(1-vz)^{-1},$$
for all $z\in\mathbb{C}$ such that 
$$|z|<\min\{|u_i|^{-1}:1\leq i\leq r, \; u_i\neq 0\}\cup\{|v|^{-1}\}.$$ 
Multiplying out to clear fractions, this is to say that
$$\sum_{1\leq i\leq r}(1-vz)\prod_{j\neq i}(1-u_jz)=r\prod_{1\leq j\leq r}(1-u_jz).$$
Equating coefficients of $z^k$ in this polynomial of degree $r$
 in $z$ yields
$$\sum_{1\leq i\leq r}e_k(u_1,\ldots, u_{i-1},v,u_{i+1},\ldots,
u_r)=re_k(u_1,\ldots, u_r),$$
where $e_k$ is the $k$th elementary symmetric function. Cancelling
terms just involving the $u_i$ (which comprise altogether $r-k$ copies of
$e_k(u_1,\ldots, u_r)$ on the left-hand side of the equation:
each $k$-subset of $\{u_1,\ldots, u_r\}$ occurs in $r-k$ of the sets
$\{u_1,\ldots, u_{i-1},v,u_{i+1},\ldots, u_r\}$) and factoring
out the resulting common factor of $v$ (on the left-hand side of the equation, each $(k-1)$-subset of $\{u_1,\ldots, u_r\}$
occurs in $r-k+1$ terms),
$$(r-k+1)ve_{k-1}(u_1,\ldots, u_r)=ke_k(u_1,\ldots, u_r).$$ 
By this recursive formula we obtain
$$e_k(u_1,\ldots, u_r)=\frac{r-k+1}{k}ve_{k-1}=\binom{r}{k}v^k.$$
This implies that $u_1,\ldots, u_r$ are uniquely determined as the
roots of the polynomial $(z-v)^r$, i.e., $u_i=v$ for each $1\leq i\leq r$.
\end{proof}

\begin{proof} {\it (of Theorem \ref{TG Hamming}.)}  In one direction, given that $w_{a,b}=w$ for $a\neq b$
  and $w_{a,a}=y$, the evaluation of the Tutte polynomial follows from that of the
  monochromial given in equation \eqref{eqn: mono poly}
  in Example \ref{ex: monochrome} above with $x=y/w$.

In the other direction, suppose that there are constants
$\alpha,\beta,\gamma, x,y$ such that $F(G;\mathbf{w})=F(G)$ satisfies
the relations \eqref{TG} for a Tutte-Grothendieck invariant.
By checking that this is indeed the case for the three families of
graphs $X_m, Y_m, Z_m$ ($m\in\mathbb{N}$) itemized below 
the desired conclusion is
reached. Each of these families of graphs 
possess the following virtues:  (i) the graphs obtained by contracting
or deleting an edge are
of the same form or belong to one of the other families, and (ii) it is
possible to write down $F$ as given by the partition function
\eqref{eqn: S} as a succinct formula, thereupon to substitute this
formula into the contraction-deletion recurrence of \eqref{TG}, and
finally to avail oneself of Lemma
\ref{lemma: mth powers} (or a variant of this lemma). 

For $m\in\mathbb{N}$ consider: 
\begin{itemize}

\item[(i)] $Y_m$, the graph on one vertex with $m$ loops.
  $F(Y_1)=\sum_{a\in Q}w_{a,a}=qy$. The relation
  $F(Y_m)=yF(Y_{m-1})=qy^m$ is used to show that
  $w_{a,a}=y$ for each $a\in Q$.

\item[(ii a)] $X_m$, the graph on two vertices connected by $m$ parallel
  edges. $F(X_1)=\sum_{a,b\in Q}w_{a,b}=qx$. The relation
  $F(X_m)=\alpha F(Y_{m-1})+\beta F(X_{m-1})$ for $m\geq 2$ is used to show that $\{w_{a,b}:a,b\in Q\}$ contains
  at most two distinct values $y,w$: there is $S$ with
  $\{(a,a):a\in Q\}\subseteq S\subseteq Q\times Q$ such that $w_{a,b}=y$ for $(a,b)\in S$, and $w_{a,b}=w$ otherwise.

\item[(ii b)] $X_m^n$,
the graph $X_m$ with $n$ edges oriented in one direction, $m-n$ in the
other. That $F(X_m)$ is independent of any orientation of the edges of $X_m$
(giving a graph $X_m^n$) is used to show that $w_{a,b}=w_{b,a}$ for all $a,b\in Q$,
i.e., the set $S$ defined in (ii a) is closed under the involution $(a,b)\mapsto (b,a)$.

\item[(iii)] $Z_m$, the star graph with $m$ edges (one vertex degree $m$, and $m$
vertices degree $1$). The relation $F(Z_m)=xF(Z_{m-1})=qx^m$ is used to show that $\#\{b\in Q:(a,b)\in S\}$ is
independent of $a$, whereby it follows from (ii b) that either
$S=\{(a,a):a\in Q\}$ or $S=Q\times Q$.

\end{itemize}
\end{proof}

We now know how to recognize a Tutte-Grothendieck invariant. To aid
our search amongst graph polynomials of the sort defined in Section
\ref{sec: Petersen graph poly} we shall use instruments from Fourier
analysis, a subject to which we now turn.

\section{Fourier analysis on finite Abelian groups} \label{sec: Fourier}

\subsection{The algebra $\mathbb{C}^Q$}

Let $Q$ be an additive Abelian group of order $q$. In later sections
$Q=\mathbb{Z}_q$, the integers under addition modulo $q$. %, sometimes $\mathbb{F}_q$.

The set $\mathbb{C}^Q$ of
functions $f:Q\rightarrow\mathbb{C}$ forms a $q$-dimensional Hermitian
inner product space. The inner product is defined by
$$\langle f,g\rangle=\sum_{a\in Q}f(a)\ol{g(a)},$$
the bar denoting complex conjugation. 
The {\em $\ell_2$-norm} is defined by $\|f\|_2=\langle
f,f\rangle^{\frac12}$ and defines a metric on the space $\mathbb{C}^Q$.

The space $\mathbb{C}^Q$ has an orthonormal basis of {\em indicator functions}  $\{\delta_a:a\in Q\}$, 
$$\delta_a(b)=\begin{cases} 1 & a=b,\\ 0 & a\neq b.\end{cases}$$
%By definition, $$f=\sum_{a\in Q}f(a)\delta_a.$$
%In other words, $f(a)=\langle f,\delta_a\rangle$.

There are several definitions of 
multiplication that make $\mathbb{C}^Q$ an algebra:
\begin{itemize}
\item[(i)] {\em Pointwise product}
$$f\cdot g(a)=f(a)g(a).$$
\item[(ii)] {\em Convolution}
$$f\ast g(a)=\sum_{b\in Q}f(a)g(b-a).$$
\item[(iii)] {\em Cross-correlation}
$$f\star g(a)=\sum_{b\in Q}\ol{f(a)}g(b+a).$$
\end{itemize}
The effect of these operations on the indicator functions is
as follows:
 $$\delta_a\cdot\delta_b=\delta_a(b)\delta_a,\hspace{1cm}\delta_a\ast\delta_b=\delta_{a+b},\hspace{1cm}\delta_a\star\delta_b=\delta_{b-a}.$$

The Abelian group $Q$ has {\em dual group} equal to the set of {\em
  characters} of $Q$ under pointwise multiplication. For each $c\in
Q$, the character
$\chi_c: Q\rightarrow\mathbb{C}^\times$ is a group homomorphism:
$\chi_c(a+b)=\chi_c(a)\chi_c(b)$ for
all $a,b\in Q$. 
The multiplicative group of characters of $Q$ is isomorphic to the
additive group $Q$.
(This is only true when $Q$ is a finite Abelian group, and, for
 the applications later in this article, is the reason why only finite Abelian groups are considered.)

The set $\{\chi_c:c\in Q\}$ forms an orthogonal basis for
$\mathbb{C}^Q$, with $\langle\chi_a,\chi_b\rangle=q\delta_a(b)$.

In the algebra $\mathbb{C}^Q$, 
$$\chi_a\cdot\chi_b=\chi_{a+b},\hspace{1cm}\chi_a\ast\chi_b=q\delta_a(b)\chi_a=\chi_{a}\star\chi_b.$$

Supposing the additive group $Q$ has the further structure of a ring
(such as $\mathbb{Z}_q$ with addition and multiplication modulo $q$), 
a {\em generating character} $\chi$ satisfies $\chi_a(b)=\chi(ab)$ for
all $a,b\in Q$.
When $Q=\mathbb{Z}_q$, the character $\chi$ defined by
$\chi(a)=e^{2\pi ia/q}$ (or $e^{2\pi i ca/q}$ for any fixed $c$
coprime with $q$) is a generating character.

\subsection{The Fourier transform}

The evaluation of the Fourier transform of a function at a point is
the projection of the function onto a character: 
$$\wh{f}(b)=\langle f,\chi_b\rangle=\sum_{a\in Q}f(a)\chi_b(-a),$$
i.e., $$\wh{f}=\sum_{b\in Q}f(b)\chi_{-b}.$$
Orthogonality of the basis $\{\chi_c:c\in Q\}$ yields:

\begin{itemize} 

\item[(i)] the Fourier inversion formula, 
%$\wh{\wh{f}}(b)=\langle\wh{f},\chi_{b}\rangle=qf(-b)$, i.e., 
$$f(a)=q^{-1}\langle\wh{f},\chi_{-a}\rangle=q^{-1}\sum_{b\in Q}\wh{f}(b)\chi_b(a),$$
i.e., the Fourier transform may be regarded as a change of basis from
indicators to characters:  
$$f=\sum_{a\in Q} f(a)\delta_a=q^{-1}\sum_{b\in Q}\wh{f}(b)\chi_b.$$

\item[(ii)] Plancherel's formula,
 $$\langle \wh{f},\wh{g}\rangle=q\langle f,g\rangle.$$

\item[(iii)] Parseval's formula,
 $$\|f\|_2^2=q^{-1}\|\wh{f}\|_2^2.$$
Thus the normalized Fourier transform $f\mapsto q^{-\frac12}\wh{f}$ is
  a unitary transformation, giving an isometry of the metric space
  $\mathbb{C}^Q$.
\end{itemize}

The Fourier transform is an isomorphism of the algebra
$(\mathbb{C}^Q,\ast)$ with the algebra $(\mathbb{C}^Q,\cdot)$:
\be\label{eqn: point to conv}\wh{f\cdot g}=q^{-1}\wh{f}\ast\wh{g},\hspace{1cm}\wh{f\ast g}=\wh{f}\cdot\wh{g},\hspace{1cm}\wh{f\star g}=\ol{\wh{f}}\cdot\wh{g},\ee
and in particular 
 $$\wh{f\star f}=|\wh{f}|^2$$
(the finite version of the Wiener-Khintchine formula). 
That the Fourier transform is an isometry carrying convolution to pointwise
multiplication makes it useful in the analysis of random walks on
Cayley graphs on $Q$, where steps on the graph correspond to addition
of group elements -- see for example \cite{Terras99} and
the references therein. 
To prove the formulae in \eqref{eqn: point to conv} it suffices to
determine the effect of the Fourier transform on basis
functions and then appeal to linearity and distributivity. For
example, $\wh{\delta_a\star\delta_b}=\wh{\delta_{b-a}}=\chi_{a-b}=\ol{\wh{\delta_a}}\cdot\wh{\delta_b}.$\smallskip

For an additive subgroup $P$ of $Q$, the {\em annihilator} of $P$ is defined by
$$P^\sharp=\{b\in Q:\forall_{a\in P}\; \chi_b(a)=1\}.$$
and is isomorphic to the quotient group $Q/P$.  
%By
%definition, if $b\in P^\sharp$ 
% then $\sum_{a\in P}\chi_b(a)=|P|$. Otherwise there is $c\in P$ with $\chi_b(c)\neq 1$ and, since $P$ is closed under addition, 
%$$\sum_{a\in P}\chi_b(a)=\sum_{a\in P}\chi_b(a+c)=\chi_b(c)\sum_{a\in
%   P}\chi_b(a),$$ implying that $\wh{\delta_P}(b)=0$ when $b\not\in
% P^\sharp$. 

Extend the indicator function notation from elements to subsets
$P\subseteq Q$
by setting $\delta_P=\sum_{a\in P}\delta_a$.

For our purposes, a key property of the Fourier transform is
that
 $$\wh{\delta_P}=|P|\delta_{P^\sharp}.$$

By Fourier inversion,
$$\delta_P\star f(b)=q^{-1}\langle\wh{\delta_P}\cdot\wh{f},\chi_{-b}\rangle,$$ giving the
{\em Poisson summation formula}
 $$\sum_{a\in P}f(a+b)=|P^\sharp|^{-1}\sum_{a\in P^\sharp}\wh{f}(a)\chi_b(a).$$
%We shall extend the domain of functions from elements of $Q$ to
%subsets by setting $f(P)=\sum_{a\in P}f(a)$. In this notation, the
%Poisson summation formula is
% $$f(P+b)=\frac{1}{|P^\sharp|}\wh{f}\cdot\chi_b(P^\sharp).$$ 

\subsection{The algebra $\mathbb{C}^{Q^n}$ and the polynomial ring
  $\mathbb{C}[\mathbf{x}]/(x_i^q-1:1\leq i\leq n)$}
In this section we assume that the Abelian group $Q$ 
has the further structure of a commutative ring. Let $\mathbf{Q}=Q^n$ denote the $n$-fold direct product of $Q$, which is an
Abelian group of order $q^n$ and a module over $Q$. Put a ring structure on $Q^n$ by defining
componentwise multiplication of  $\mathbf{a}=(a_1,\ldots, a_n)$, $\mathbf{b}=(b_1,\ldots, b_n)\in Q^n$,  
$$\mathbf{a}\mathbf{b}=(a_1b_1,\ldots, a_nb_n).$$
The Hermitian inner product space $\mathbb{C}^{Q^n}$ is the $n$-fold
tensor product of $\mathbb{C}^Q$: given $f_1,\ldots,
f_n\in\mathbb{C}^Q$ define
 $$f_1\otimes\cdots \otimes f_n(a_1,\ldots, a_n)=f_1(a_1)\cdots f_n(a_n),$$
and in particular $f^{\otimes n}(\mathbf{a})=f(a_1)\cdots f(a_n).$

The characters of $Q^n$ are the functions defined by
$\chi_{\mathbf{a}}=\chi_{a_1}\otimes \cdots\otimes\chi_{a_n}.$

Define the Euclidean (dot) product by
$$\mathbf{a}\cdot\mathbf{b}=a_1b_1+\cdots +a_nb_n.$$
If $\chi$ a generating character for $Q$, then 
$\chi^{\otimes n}$ is a generating character for $Q^n$:
$$\chi_{\mathbf{a}}(\mathbf{b})=\chi^{\otimes n}(\mathbf{a}\mathbf{b})=\chi(a_1b_1)\cdots\chi(a_nb_n)=\chi(\mathbf{a}\cdot\mathbf{b}).$$
 Given that $Q$ has a generating character, 
 for a submodule $\mathbf{P}$ of $Q^n$ the annihilator
$$\mathbf{P}^\sharp=\{\mathbf{b}\in Q^n:\forall_{\mathbf{a}\in \mathbf{P}}\hspace{0.3cm}\chi_{\mathbf{a}}(\mathbf{b})=1\}$$
is equal to the {\em orthogonal} submodule
$$\mathbf{P}^\perp=\{\mathbf{b}\in Q^n:\forall_{\mathbf{a}\in\mathbf{P}}\hspace{0.3cm}\mathbf{a}\cdot\mathbf{b}=0\}.$$

The Fourier transform on $Q^n$ is given by
$$\wh{f_1\otimes\cdots\otimes f_n}=\wh{f_1}\otimes\cdots\otimes\wh{f_n},$$
and in particular $\wh{f^{\otimes n}}=\wh{f}^{\otimes n}.$

It may be helpful to spell 
out the relationship between polynomials in the ring
$\mathbb{C}[\mathbf{x}]/(x_i^q-1:1\leq i\leq n)$ (where
$\mathbf{x}=(x_i:1\leq i\leq n)$ is an $n$-tuple of commuting
indeterminates) and functions in the space
$\mathbb{C}^{\mathbb{Z}_q^n}$. The aim of course is to translate
statements about the reduced graph polynomial $F^{(q)}(G;\mathbf{x})$, which
belongs to $\mathbb{C}[\mathbf{x}]/(x_v^q-1:v\in V)$, into statements
about functions in $\mathbb{C}^{\mathbb{Z}_q^V}$. The latter space has now
the advantage of familiarity and the accoutrements of a succinct
notation.

Take
$Q=\mathbb{Z}_q$, which has generating character $\chi(a)=\zeta^{a}$ for
 primitive $q$th root of unity $\zeta$.

The algebra $\mathbb{C}^{\mathbb{Z}_q^n}$ is isomorphic to
$\mathbb{C}[\mathbf{x}]/(x_i^q-1:1\leq i\leq n)$ and the
following correspondences obtain:
$$\delta_{\mathbf{a}}=\delta_{a_1}\otimes\cdots\otimes
\delta_{a_n}\hspace{1cm}\mbox{with}\hspace{1cm}\mathbf{x}^{\mathbf{a}}=\prod_{1\leq
  i\leq n}x_i^{a_i},$$
$$\chi_{\mathbf{a}}=\chi_{a_1}\otimes\cdots\otimes\chi_{a_n}\hspace{1cm}\mbox{with}\hspace{1cm}\prod_{1\leq
  i\leq n}\frac{x_i^q-1}{\zeta^{-a_i}x_i-1},$$
$$\mathbf{f}=\sum_{\mathbf{a}\in\mathbb{Z}_q^n}\mathbf{f}(\mathbf{a})\delta_{\mathbf{a}}\hspace{1cm}\mbox{with}\hspace{1cm}F(\mathbf{x})=\sum_{\mathbf{a}\in
  \mathbb{Z}_q^n}\mathbf{f}(\mathbf{a})\mathbf{x}^{\mathbf{a}}.$$ 
Finally, $$F(\zeta^{a_1},\ldots,
\zeta^{a_n})=\wh{\mathbf{f}}(\mathbf{a}),$$
and Lagrange interpolation on points $\{(\zeta^{a_1},\ldots,
\zeta^{a_n}):(a_1,\ldots,a_n)\in\mathbb{Z}_q^n\}$ is the Fourier basis change:
$$\sum_{\mathbf{a}\in \mathbb{Z}_q^n}\mathbf{f}(\mathbf{a})\mathbf{x}^{\mathbf{a}}=q^{-n}\sum_{\mathbf{a}\in\mathbb{Z}_q^n}\wh{f}(\mathbf{a})\prod_{i=1}^n\frac{x_i^q-1}{\zeta^{-a_i}x_i-1}.$$

\subsection{Weight enumerators and the Tutte polynomial}\label{sec: we}

We finish this section on Fourier analysis with a discussion of the
Tutte polynomial as a weight enumerator that gives us the opportunity
at the same time to define {\em
  flows} and {\em tensions} of a graph, which definitions are needed
for the next section.

It will be convenient to extend the domain of a function $\mathbf{f}$
on elements $\mathbf{a}\in Q^n$ to subsets $\mathbf{P}\subseteq Q^n$, setting
$$\mathbf{f}(\mathbf{P})=\sum_{\mathbf{a}\in\mathbf{P}}\mathbf{f}(\mathbf{a}).$$
The {\em Hamming weight}  of $\mathbf{a}=(a_1,\ldots, a_n)$ is
$|\mathbf{a}|=\#\{i:a_i\neq 0\}.$
The {\em Hamming weight enumerator} of $\mathbf{P}$ is defined to be the
the generating function for vectors in $\mathbf{P}$ counted according to their number of zero entries:
$$\sum_{\mathbf{a}\in\mathbf{P}}x^{n-|\mathbf{a}|}=(x\delta_0+\delta_{Q\setminus
  0})^{\otimes n}(\mathbf{P}).$$ 
The {\em complete weight enumerator} of $\mathbf{P}$ keeps account of the
number of entries equal to a given element of $Q$:
$$\sum_{\mathbf{a}\in\mathbf{P}}\prod_{a\in Q}x_a^{\#\{1\leq i\leq
  n:a_i=a\}}=\big(\sum_{a\in Q}x_a\delta_a\big)^{\otimes E}(\mathbf{P}).$$
The Hamming weight enumerator is the specialization $x_0=x$ and
$x_a=1$ for $0\neq a\in Q$.

For a submodule $\mathbf{P}$ of $Q^n$,
$$\wh{\delta_{\mathbf{P}}}=|\mathbf{P}|\delta_{\mathbf{P}^\perp}.$$
The Poisson summation formula 
$$\mathit{\mathbf{f}}(\mathbf{P}+\mathbf{b})=\frac{1}{|\mathbf{P}^\perp|}\wh{\mathbf{f}}\cdot\chi_{\mathbf{b}}(\mathbf{P}^\perp),$$
with $\mathbf{b}=\mathbf{0}$ and $\mathbf{f}=f^{\otimes n}$ gives the
duality formula between the complete weight enumerator   
of $\mathbf{P}$ (with $x_a=f(a)$) and the complete weight enumerator
of $\mathbf{P}^\perp$ (with $x_a=\wh{f}(a)$). When
$f=x\delta_0+\delta_{Q\setminus 0}$ it yields the MacWilliams duality formula for Hamming weight enumerators.

Recall that the graph $G=(V,E)$ has a fixed orientation of its edges,
with directed edge set denoted by $\overrightarrow{E}$. Represent this
ground orientation as a matrix $\gamma$ indexed by $V\times E$,
setting
$$\gamma_{v,e}=\begin{cases} +1 & \mbox{\rm if $e=(u,v)$ in $\overrightarrow{E}$,}\\
  -1 & \mbox{\rm if $e=(v,u)$ in $\overrightarrow{E}$,}\\
0 &  \mbox{\rm if $e$ is not incident with $v$.}\end{cases}$$
A {\em $Q$-tension} of $G$ is a vector $\mathbf{a}\in Q^E$ comprising the differences between endpoints
in a vertex colouring $\mathbf{c}\in Q^V$, i.e., if
$e=(u,v)$ then the $Q$-tension $\mathbf{a}$ associated with
$\mathbf{c}$ is defined by
$$a_e=\sum_{v\in V}\gamma_{v,e}c_v=c_v-c_u.$$
A {\em $Q$-flow} of $G$ is a vector $\mathbf{b}\in Q^E$ such that, for
each vertex $v$,
$$\sum_{e\in E}\gamma_{v,e}b_e=0.$$

When $G$ is planar, $Q$-flows of $G$ correspond to $Q$-tensions of the planar dual graph $G^*$. In particular, when $Q=\mathbb{F}_2$, the $\mathbb{F}_2$-flows of $G$ (cycles/ Eulerian subgraphs) correspond to $\mathbb{F}_2$-tensions (cutsets) of $G^*$. 
 
A {\em nowhere-zero} $Q$-tension is one that takes only non-zero values, and arises from a proper vertex
$Q$-colouring; similarly, a nowhere-zero $Q$-flow is a flow takes
non-zero values only
(and for plane graphs corresponds to a proper face $Q$-colouring of
the embedded graph).

If $\mathbf{P}$ is the set of $Q$-tensions of $G$ (of which there are
$q^{r(G)}$) then $\mathbf{P}^\perp$ is the set of
$Q$-flows of $G$ (of which there are $q^{n(G)}$).
With this notation, the monochromial is given by
$$\sum_{\mathbf{c}\in
  Q^V}y^{\#\{(u,v)\in\overrightarrow{E}:c_u=c_v\}}=q^{k(G)}\sum_{\mathbf{a}\in\mathbf{P}}y^{|E|-|\mathbf{a}|},$$
since there are $q^{k(G)}$ vertex $Q$-colourings yielding any given $Q$-tension. Consequently, by Example \ref{ex: monochrome}, the Hamming weight enumerator of the set $\mathbf{P}$
of $Q$-tensions of $G$ is a specialization of the Tutte polynomial: 
$$\sum_{\mathbf{a}\in\mathbf{P}}y^{|E|-|\mathbf{a}|}=(y-1)^{r(G)}T(G;\frac{y-1+q}{y-1},y).$$

By the Poisson summation formula (MacWilliams duality),
$$(y\delta_0+\delta_{Q\setminus 0})^{\otimes E}(\mathbf{P})=q^{-n(G)}[(y-\!1+\!q)\delta_0+(y-\!1)\delta_{Q\setminus
    0}]^{\otimes E}(\mathbf{P}^\perp).$$
Putting $x=\frac{y-\!1+\!q}{y-\!1}$, the Hamming weight enumerator of
the set $\mathbf{P}^\perp$ of $Q$-flows of $G$ is given by
 $$\sum_{\mathbf{b}\in\mathbf{P}^\perp}x^{|E|-|\mathbf{b}|}=(x-1)^{n(G)}T(G;x,\frac{x-1+q}{x-1}).$$

A corollary of Theorem \ref{TG Hamming} is that if an evaluation of the complete weight
enumerator of $Q$-tensions (or $Q$-flows) is a Tutte-Grothendieck
invariant (an evaluation of the
Tutte polynomial with a certain simple type of prefactor) 
then it is in fact an evaluation of the Hamming weight
enumerator.
In fact, the proof of Theorem \ref{TG Hamming} says the same is true
of any class of graphs that contains multiple loops on one vertex,
multiple parallel
edges between two vertices, and stars whose central vertex is of
arbitrary degree. This notably includes the
class of planar graphs.

There are nevertheless (infinite) classes of graphs for which an
evaluation of the complete weight enumerator of $Q$-tensions of $G$
coincides with the value of a Tutte-Grothendieck invariant and yet is not an evaluation of the
Hamming weight enumerator.

For example, if $G=(V,E)$ is the line graph of a plane cubic
graph then a result ultimately due to Penrose \cite{P71} (but see \cite{EG96} for
a full account) is that 
$$\sum_{\mathbf{c}\in\mathbb{Z}_3^V}0^{\#\{(u,v)\in\overrightarrow{E}:
  c_u=c_v\}}(-1)^{\#\{(u,v)\in
  \overrightarrow{E}:c_v-c_u=-1\}}=(-1)^{|V|}P(G;3),$$
i.e., the complete weight enumerator of $\mathbb{Z}_3$-tensions of $G$
with $x_0=0, x_1=1, x_{-1}=-1$ is an evaluation of the Tutte
polynomial. However, since the class of line graphs of plane cubic
graphs is not
closed under deletion or contraction, one is prevented from calling this a
Tutte-Grothendieck invariant.

\section{Polynomials akin to the graph polynomial}\label{sec: graph polynomials}

 Suppose $F^{(q)}(G;\mathbf{x})\in\mathbb{C}[\mathbf{x}]/(x_v^q-1:v\in V)$ is a graph polynomial of the general
 form
\begin{align*}F^{(q)}(G;\mathbf{x}) & = \prod_{(u,v)\in\overrightarrow{E}}\sum_{(a,b)\in
  \mathbb{Z}_q^2}f(a,b)x_u^ax_v^b\\
 & =\sum_{\mathbf{c}\in (\mathbb{Z}_q^2)^E}f^{\otimes
  E}(\mathbf{c})\prod_{(u,v)\in \overrightarrow{E}}x_u^{c_{u,e}}x_v^{c_{v,e}},\end{align*}
where $\mathbf{c}=(c_e:e\in E)$, $c_e=(c_{u,e},c_{v,e})$ for edge $e$
directed as $(u,v)$ in $\overrightarrow{E}$, and
$f^{\otimes E}(\mathbf{c})=\bigotimes_{e\in E}f(c_{u,e},c_{v,e}).$
The graph polynomial of Petersen et al. introduced in Section
\ref{sec: Petersen graph poly} is the case $f(1,0)=1,
f(0,1)=-1$ and $f(a,b)=0$ otherwise. (Henceforth the name ``Petersen's
   graph polynomial'' will be used when it needs to be distinguished.)

In this section we address the following questions: 

\begin{itemize}

\item[(A)] When is the partition function  of the vertex colouring (states) model\,\footnote{ The vertex colouring model assigns weight $F^{(q)}(G;(\zeta^{d_v}:v\in V)\,)$ to a given vertex colouring $\mathbf{d}\in\mathbb{Z}_q^V$. In terms of graph homomorphisms, this vertex colouring model corresponds to considering $\mathbf{d}$ as a homomorphism from $G$ to a weighted directed graph $H$ on vertex set $\mathbb{Z}_q$, with an edge $(c,d)$ having weight 
$$\sum_{a,b}f(a,b)\zeta^{ac+bd}=\wh{f}(c,d).$$
The total weight of the homomorphism $\mathbf{d}: G\rightarrow H$ is
the product of all the weights on $(d_u,d_v)$ for edges $(u,v)$ of
$G$, i.e., $\wh{f}^{\otimes E}(\mathbf{c})$ where
$\mathbf{c}\in(\mathbb{Z}_q^2)^E$ is defined by $(c_{u,e},
c_{v,e})=(d_u,d_v)$. The partition function in question~(A) is a sum over all
homomorphisms [$\mathbf{d}$, encoded by $\mathbf{c}\in(\mathbb{Z}_q^2)^E$] weighted in this way.}
$$\sum_{\mathbf{d}\in \mathbb{Z}_q^V}F^{(q)}(G;(\zeta^{d_v}:v\in V)\,)=q^{|V|}[\mathbf{x}^{\mathbf{0}}]F^{(q)}(G;\mathbf{x})$$
a Tutte-Grothendieck invariant (an evaluation of the Tutte polynomial)? 

\item[(B)] When is the squared
$\ell_2$-norm 
$$\|F^{(q)}(G;\mathbf{x})\|_2^2=\sum_{\mathbf{a}\in\mathbb{Z}_q^V}\left|[\mathbf{x}^{\mathbf{a}}]F^{(q)}(G;\mathbf{x})\right|^2$$
a Tutte-Grothendieck invariant? 

\item[(C)] What are the equivalents of Theorems \ref{thm: Alon Tarsi} and \ref{thm: Tarsi} in this more general case? 
\end{itemize}

By Parseval's formula,
$$\|F^{(q)}(G;\mathbf{x})\|_2^2 = q^{-|V|}\sum_{\mathbf{d}\in\mathbb{Z}_q^V}|F^{(q)}(\zeta^{d_v}:v\in V)|^2,$$
where, writing $\mathbf{c}$ for the vector with entries $(c_{u,e},c_{v,e})=(d_u,d_v)$,
$$|F^{(q)}(G;\zeta^{d_v}:v\in V)|^2= |\wh{f}^{\,\otimes E}(\mathbf{c})|^2.$$ 
Since $|\wh{f}|^2=\wh{f\star f}$, this implies that the $\ell_2$-norm
of $F^{(q)}(G;\mathbf{x})$ is equal to the constant term of the polynomial
$\wt{F}^{(q)}(G;\mathbf{x})$ in $\mathbb{C}[\mathbf{x}]/(x_v^q-1:v\in V)$ defined by
$$\wt{F}^{(q)}(G;\mathbf{x})=\prod_{(u,v)\in\overrightarrow{E}}f\star
f(a,b)x_u^ax_v^b.$$

For example, the $\ell_2$-norm of Petersen's graph polynomial
$$F^{(q)}(G;\mathbf{x})=\prod_{(u,v)\in\overrightarrow{E}}(x_u-x_v)\;\;\bmod
(x_v^q-1:v\in V)$$
is the constant term of the polynomial
\begin{align*}
|F^{(q)}(G;\mathbf{x})|^2 & =\prod|x_u-x_v|^2=\prod
(x_u-x_v)(x_u^{-1}-x_v^{-1}),\\
& = \prod (2-x_ux_v^{-1}-x_u^{-1}x_v)\;\;\bmod
(x_v^q-1:v\in V).\end{align*}
[This calculation uses the correspondence of the ideal $(x_v^q-1:v\in V)$ with the
algebraic variety of points $(\zeta^{d_v}:v\in V)$, i.e.,
indeterminates $x_v$ are roots of unity, for which complex
conjugates are the same as multiplicative inverses.]

Let $M=\{(a,a):a\in Q\}$ be the submodule of $Q\times Q$ comprising
monochromatic pairs. The orthogonal submodule is $M^\perp=\{(a,-a):a\in Q\}$. 
By Theorem \ref{TG Hamming}, $[\mathbf{x}^\mathbf{0}]F^{(q)}(G;\mathbf{x})$ is a
Tutte-Grothendieck invariant if and only if there are constants $y,w$
such that $\wh{f}=y\delta_M+w\delta_{Q\times Q\setminus M}$. By the above
remarks, $\|F^{(q)}(G;\mathbf{x})\|_2^2$ is a Tutte-Grothendieck invariant
if and only if 
$$\wh{f\star f}=y\delta_{M}+w\delta_{Q\times Q\setminus M}.$$
By Fourier inversion, 
this is the case if and only if $f\star f=(y+(q-1)w)\delta_0+(y-w)\delta_{M^\perp\setminus 0}$. 

\begin{proposition} \label{prop: constant} The
constant term of $F^{(q)}(G;\mathbf{x})$ is a Tutte-Grothendieck
invariant if and only if
$$F^{(q)}(G;\mathbf{x})=\prod_{(u,v)\in\overrightarrow{E}}\big[
  y+(q-1)w+(y-w)(x_u^{q-1}x_v+\cdots + x_ux_v^{q-1})\big]\bmod(x_v^q-1:v\in
V),$$
in which case
$$[\mathbf{x}^\mathbf{0}]F^{(q)}(G;\mathbf{x})=(qw)^{n(G)}(y-w)^{r(G)}T(G;\frac{y-(q-1)w}{y-w},\frac{y}{w}).$$
\end{proposition}
For example, when $y=0, w=1$ and $q=3$ this says that 
$$\prod_{(u,v)\in\overrightarrow{E}}(2-x_ux_v^{2}-x_u^{2}x_v)\bmod(x_v^3-1:v\in V)$$
has constant term $3^{n(G)}(-1)^{r(G)}T(G;-2,0)=3^{|E|-|V|}P(G;3)$.

That the constant term of the polynomial defined in Proposition \ref{prop: constant} is a Tutte
polynomial evaluation can be seen by inspection since,
for $a,b\in\mathbb{Z}_q$,
$$\zeta^{(q-1)a+b}+\zeta^{(q-2)a+2b}+\cdots +\zeta^{a+(q-1)b}=\begin{cases} -1 & a\neq b\\
 q-1 & a=b,\end{cases}$$
so that in this case
\begin{align*}
q^{|V|}[\mathbf{x}^{\mathbf{0}}]F^{(q)}(G;\mathbf{x}) &
  =\sum_{\mathbf{a}\in\mathbb{Z}_q^V}F^{(q)}(G;(\zeta^{a_v}:v\in
  V)\,)\\
& =
  \sum_{\mathbf{c}\in\mathbb{Z}_q^V}(qy)^{\#\{(u,v)\in\overrightarrow{E}:c_u=c_v\}}(qw)^{\#\{(u,v)\in\overrightarrow{E}:c_u\neq
    c_v\}}.\end{align*}

Whereas Proposition \ref{prop: constant} limits the graph polynomials which have
constant term equal to an evaluation of the Tutte polynomial to a
single family -- giving a rather dull answer to question~(A) above -- the possible choices for $f$ defining
$F^{(q)}(G;\mathbf{x})$ so that the $\ell_2$-norm $\|F^{(q)}(G;\mathbf{x})\|_2^2$ is a
Tutte-Grothendieck invariant are unlimited -- making the answer to
question~(B) potentially equally as dull.
The criterion $|\wh{f}|^2=y\delta_M+w\delta_{Q\times Q\setminus M}$
[or $f\star f=(y+(q-1)w)\delta_0+(y-w)\delta_{M^\perp\setminus 0}$]
can be satisfied by taking $\wh{f}=\sum_{a\in  Q}z_{a,b}\delta_{(a,b)}$ for any complex numbers $z_{a,b}$ that
satisfy $|z_{a,a}|^2=y$ if $a=b$ and $|z_{a,b}|^2=w$ otherwise.

Nonetheless, it seems worth describing a family of polynomials which
contains Petersen's graph polynomial as a special case and in some sense naturally
generalizes it.  In this family it is also possible to give a
meaningful answer to question (C) asking for equivalents to
Theorems \ref{thm: Alon Tarsi} and \ref{thm: Tarsi}.

\subsection{A family of polynomials containing the graph polynomial}

 Suppose that supp$(f)\subseteq\{(a,b):a+sb=t\}$ for some constants
$s,t\in\mathbb{Z}_q$. Then
\begin{align*}
F^{(q)}(G;(\zeta^{d_v}:v\in
V)\,) & =\prod_{(u,v)\in\overrightarrow{E}}\sum_{(a,b)\in
  \mathbb{Z}_q^2}f(a,b)\zeta^{ad_u+bd_v}\\
 & =\sum_{\mathbf{c}\in(\mathbb{Z}_q^2)^E}f^{\otimes
  E}(\mathbf{c})\prod_{(u,v)\in\overrightarrow{E}}\zeta^{c_{u,e}d_u+c_{v,e}d_v}.\end{align*}
The equation $f(a,b)=f(t-sb,b)=:g(b)$ defines $g\in\mathbb{C}^{\mathbb{Z}_q}$ and
 the sum over
$\mathbf{c}\in(\mathbb{Z}_q^2)^E$ can be rewritten as a sum over $\mathbf{b}\in\mathbb{Z}_q^E$. 
In particular,  $s=1$ when the polynomial $\sum_{a,b}f(a,b)x_u^ax_v^b$ is homogeneous.

Given that $a_e+sb_e=t$, we have
$a_{e}d_u+b_ed_v=(t-sb_e)d_u+b_ed_v=b_e(d_v-sd_u)+td_u.$
For $e=(u,v)\in\overrightarrow{E}$, define $S:\mathbb{Z}_q^V\rightarrow\mathbb{Z}_q^E$  by
$$(S\mathbf{d})_e=d_v-sd_u$$ and $T:\mathbb{Z}_q^V\rightarrow\mathbb{Z}_q^E$ by $$(T\mathbf{d})_e=td_u.$$

For $\mathbf{b}\in\mathbb{Z}_q^E$, the transpose $S^\top$ is given by
$$(S^\top\mathbf{b})_v=\sum_{e=(u,v)\in\overrightarrow{E}}b_e\;-\hspace{0.5cm}s\sum_{e=(v,u)\in\overrightarrow{E}}b_e$$ 
and $$(T^\top\mathbf{b})_v=t\sum_{e=(v,u)\in\overrightarrow{E}}b_e.$$
When $s=1$ (which is the case for Petersen's graph
polynomial) the linear transformation $S$ is the {\em coboundary} and $S^\top$ the
{\em boundary}. Here the submodule $\ker(S^\top)$ comprises the
 $\mathbb{Z}_q$-flows of $G$ and im$(S)$ the  $\mathbb{Z}_q$-tensions
of $G$.

We have
\begin{align*}F^{(q)}(G;\mathbf{x}) & =
  \prod_{(u,v)\in\overrightarrow{E}}\;\;\sum_{b\in\mathbb{Z}_q}g(b)x_u^{t-sb}x_v^b\\
 & =\sum_{\mathbf{b}\in\mathbb{Z}_q^E}\prod_{e=(u,v)\in \overrightarrow{E}}g(b_e)x_u^{t-sb_e}x_v^{b_e}\\
& = \sum_{\mathbf{b}\in\mathbb{Z}_q^E}g^{\otimes
    E}(\mathbf{b})\prod_{v\in
    V}x_v^{S^\top\mathbf{b}+T^\top\mathbf{1}},\end{align*}
where $\mathbf{1}$ is the all-one vector in
$\mathbb{Z}_q^V$. ($T^\top\mathbf{1}$ is $t$ times the outdegree score
of $\overrightarrow{E}$.) \smallskip

The following theorem provides an answer to the question~(A) posed in the
previous section, and more.

\begin{theorem}
If $S^\top\mathbf{b}=\mathbf{a}-T^\top\mathbf{1}$ then 
 $$[\mathbf{x}^{\mathbf{a}}]F^{(q)}(G;\mathbf{x}) = g^{\otimes
  E}(\ker(S^\top)+\mathbf{b}),$$
a complete coset weight enumerator of $\ker(S^\top)$.

In particular,  the coefficient
$[\mathbf{x}^{T^\top\mathbf{1}}]F^{(q)}(G;\mathbf{x})$ is an evaluation of the
 complete weight
  enumerator of $\ker(S^\top)$ (and of im$(S)$). 
\end{theorem} 
For example, in Petersen's graph polynomial, where $g=\delta_0-\delta_1$, 
$$[\mathbf{x}^{T^\top\mathbf{1}}]\prod_{(u,v)\in\overrightarrow{E}}(x_u-x_v)\bmod(x_v^q-1:v\in
V) = \sum_{(q,1)\mbox{\rm \small -flows } \, \mathbf{b}}0^{\#\{e\in
  E:b_e=-1\}}(-1)^{\#\{e\in E:b_e=1\}},$$
where a $(q,1)$-flow is a $\mathbb{Z}_q$-flow taking values only in
$\{0,1,-1\}$ (and here the sum need only range over those taking
values in $\{0,1\}$).

When $s=1$ (for which $S$ is the coboundary, ${\rm im}(S)$ the set of
$\mathbb{Z}_q$-tensions, $\ker(S^\top)$ the set of
$\mathbb{Z}_q$-flows) and
$$F^{(q)}(G;\mathbf{x})=\prod_{(u,v)\in\overrightarrow{E}}\;\sum_{b\in\mathbb{Z}_q}g(b)x_u^{t-b}x_v^b,$$
 the coefficient
$[\mathbf{x}^{T^\top\mathbf{1}}]F^{(q)}(G;\mathbf{x})$ is an evaluation of the Tutte
polynomial if and only if $g=x\delta_0+\delta_{\mathbb{Z}_q\setminus
  0}$ (by Theorem \ref{TG Hamming}; this is the case covered by
Proposition \ref{prop: constant}). If $g$
does not take this form then the coefficient
$[\mathbf{x}^{T^\top\mathbf{1}}]\ol{F}(\mathbf{x})$ is not an
evaluation of the Hamming
weight enumerator of $\mathbb{Z}_q$-flows but of some other
specialization of the complete weight enumerator. \medskip

To find the $\ell_2$-norm, observe that, for $\mathbf{d}\in\mathbb{Z}_q^V$,
\begin{align*}
  F^{(q)}(G;(\zeta^{d_v}:v\in
V)\,) & = \sum_{\mathbf{b}\in\mathbb{Z}_q^E}g^{\otimes
  E}(\mathbf{b})\zeta^{(S^\top\mathbf{b})\cdot \mathbf{d}+T^\top\mathbf{1}\cdot\mathbf{d}}\\
 & = \sum_{\mathbf{b}\in\mathbb{Z}_q^E}g^{\otimes
  E}(\mathbf{b})\zeta^{\mathbf{b}\cdot S\mathbf{d}+\mathbf{1}\cdot
  T\mathbf{d}}\\
& =\zeta^{\mathbf{1}\cdot T\mathbf{d}}\wh{g}^{\otimes E}(-S\mathbf{d}),\end{align*}
and
$$|F^{(q)}(G;(\zeta^{d_v}:v\in V)\,)|^2=|\wh{g}^{\otimes{E}}(-S\mathbf{d})|^2=|\wh{g}^{\otimes{E}}(S\mathbf{d})|^2.$$

 By Parseval's formula,
\begin{align*}\|F^{(q)}(G;\mathbf{x})\|_2^2  & = q^{-|V|}\sum_{\mathbf{d}\in\mathbb{Z}_q^V}|\wh{g}^{\otimes
  E}(S\mathbf{d})|^2.\\
& 
=q^{-|V|}|\ker(S)|\sum_{\mathbf{b}\in\mbox{\rm \small
    im}(S)}(|\wh{g}|^2)^{\otimes E}(\mathbf{b}).\end{align*}

By the Poisson summation formula, and
using  ${\rm im}(S)^\perp=\ker(S^\top)$,
$|\ker(S)|=q^{|V|}/|{\rm im(S)}|$, we deduce the following, which
provides an answer to question (C).
 \begin{theorem}\label{thm: l2 norm}
If $$F^{(q)}(G;\mathbf{x})=\prod_{(u,v)\in\overrightarrow{E}}\;\sum_{b\in
    \mathbb{Z}_q}g(b)x_u^{t-sb}x_v^{b},$$
then
\begin{align*}\|F^{(q)}(G;\mathbf{x})\|_2^2 & =\frac{1}{|{\rm im}(S)|}\sum_{\mathbf{b}\in{\rm
    im}(S)}|\wh{g}^{\otimes E}|^2(\mathbf{b})\\
& =\sum_{\mathbf{b}\in\ker(S^\top)}(g\star g)^{\otimes
    E}(\mathbf{b}),\end{align*}
where as usual $S:\mathbb{Z}_q^V\rightarrow\mathbb{Z}_q^E$ is defined by
$(S\mathbf{d})_e=d_v-sd_u$ for $e=(u,v)\in\overrightarrow{E}$.
\end{theorem}

\begin{example} Petersen's graph polynomial modulo $(x_v^q-1:v\in V)$ has $s=1=t$,
$g=\delta_0-\delta_1$, $g\star g=2\delta_0-\delta_1-\delta_{-1}$. The
transformation $S:\mathbb{Z}_q^V\rightarrow\mathbb{Z}_q^E$ is the
coboundary operator, $S^\top$ the boundary, $\ker(S^\top)$ the
set of $\mathbb{Z}_q$-flows of $G$. 
This gives Tarsi's result, Theorem \ref{thm: Tarsi}, that the $\ell_2$-norm of Petersen's graph
polynomial modulo $(x_v^q-1:v\in V)$ is equal to 
$$(-1)^{|E|}\sum_{\mathbf{b}\in\{-1,0,1\}^E\cap\ker(S^\top)}(-2)^{\#\{e\in E: b_e=0\}},$$
where the sum is over {\em $(q,1)$-flows} of $G$. 
\end{example}

\begin{example} The polynomial
$$\prod_{uv\in E}(x_u+x_v)$$
is a generating function for score vectors of orientations of $G$, and
as such its number of non-zero coefficients turns out to be equal to
$T(G;2,1)$, the number of forests of $G$. (See for example \cite{BO92}).
By Theorem \ref{thm: l2 norm} with  $g=\delta_0+\delta_1, g\star g=2\delta_0+\delta_{1}+\delta_{-1}$, when this polynomial is reduced modulo $(x_v^3-1:v\in V)$ it has $\ell_2$-norm
equal to $T(G;2,4)$. Determining how many
non-zero coefficients the polynomial has (its $\ell_0$-norm) when reduced modulo
$(x_v^q-1:v\in V)$ includes as a subproblem determining whether a
graph is {\em $\mathbb{Z}_q$-connected}, a notion defined in \cite{JLPT92}. 
\end{example} 

Theorem \ref{TG Hamming} applied to the result of Theorem~\ref{thm: l2
  norm} has the following consequence, answering
question~(B).
\begin{corollary}\label{cor: l2 Tutte}
The $\ell_2$-norm $\|F^{(q)}(G;\mathbf{x})\|_2^2$ of the polynomial
defined in Theorem \ref{thm: l2 norm} is an
evaluation of the Tutte polynomial $T(G;x,y)$ with $(x-1)(y-1)=q$  if
and only if $s=1$ and $g\star g$, equivalently $|\wh{g}|^2$, is constant on
$\mathbb{Z}_q\setminus 0$.
\end{corollary}
 
We finish with three examples of functions $g$ satisfying the
conditions of Corollary \ref{cor: l2 Tutte}, yielding families of polynomials that have
$\ell_2$-norm equal to a Tutte-Grothendieck invariant. 

A {\em $(q,k,\ell)$-difference set} in an Abelian group $Q$ is a subset $P$ of size $k$
with the property that $\#\{a,b\in P:a-b=c\}=\ell$ for each $c\in
Q\setminus 0$. For example, $Q\setminus 0$ is a
$(q,q-1,q-2)$-difference set. All
non-zero $c$ have exactly $q-2$ ways of being written as $a-b$ for
$a,b\in Q\setminus 0$ since for given $a\in Q\setminus\{0,c\}$ there is a
unique $b\in Q\setminus \{0,c\}$ with $a-b=c$. %A more interesting
%example of a difference set is when $q$ is a prime power congruent to
%$3$ modulo $4$ and to take the set $P$ of non-zero squares in the
%field $\mathbb{F}_q$. This gives a
%$(q,\frac{q-1}{2},\frac{q-3}{4})$-difference set since $P$ is a
%multiplicative subgroup of $\mathbb{F}_q$ and $-1\not\in P$. 

Note that a
function is constant on non-zero values if and only if the same is true
of its Fourier transform: if $f=t\delta_0+\delta_{Q\setminus 0}$ then
$\wh{f}=(t\!-\!1\!+\!q)\delta_0+(t\!-\!1)\delta_{Q\setminus 0}$. This
fact, 
together with the equation 
$\delta_P\star\delta_P=\sum_{c\in Q}\#\{a,b\in P:a-b=c\}\delta_c$,
implies that the Fourier transform $\wh{\delta_P\star\delta_P}=|\wh{\delta_P}|^2$ is constant on
$Q\setminus 0$ if and only if $P$ is a $(q,k,\ell)$-difference set in $Q$,  i.e., $\delta_P\star\delta_P=k\delta_0+\ell\delta_{Q\setminus 0}$.

 \begin{example} If $g=\delta_P$ for some $P\subseteq
\mathbb{Z}_q$, or
more generally $g=\delta_P+r\delta_{\mathbb{Z}_q\setminus P}$ for any constant
  $r$, then $|\wh{g}|^2$ is constant on $\mathbb{Z}_q\setminus
0$ if and only if $P$ is a difference set in
$\mathbb{Z}_q$. When $P=\mathbb{Z}_q\setminus 0$ this is the family of
polynomials described in Proposition \ref{prop:
  constant} whose constant terms were also Tutte-Grothendieck invariants. 
\end{example}

A {\em $(q,k,\ell,m)$-partial difference set} in $Q$ is a subset $P$ of size $k$
with the property that
$\delta_P\star\delta_P=k\delta_0+\ell\delta_{P\setminus
  0}+m\delta_{Q\setminus (P\cup 0)}$. For example, a subgroup $P$ of
size $k$ is a $(q,k,k,0)$-partial difference set.

\begin{example}  If
$P\subseteq\mathbb{Z}_q\setminus 0$ and
$g=\delta_P-\delta_{\mathbb{Z}_q\setminus (P\cup 0)}$ then  then $|\wh{g}|^2$ is constant on $\mathbb{Z}_q\setminus
0$  iff $q$ is odd and $P$ is a Paley difference set or partial
  difference set, i.e. $|P|=(q-1)/2$ and 
$$\delta_P\star\delta_P=\begin{cases} \frac{q-1}{2}\delta_0+\frac{q-5}{4}\delta_P+\frac{q-1}{4}\delta_{\mathbb{Z}_q\setminus(P\cup
  0)} \\
\frac{q-1}{2}\delta_0+\frac{q-3}{4}\delta_{\mathbb{Z}_q\setminus 0}, \end{cases}$$
according as $q\equiv \pm 1\pmod 4$. (For odd prime $q$, the set of non-zero squares in $\mathbb{Z}_q$ is an example of such a $P$.)
\end{example}

\begin{example} When $q$ is prime and
$g=\sum_{a\in\mathbb{Z}_q}\psi(a)\delta_a$ for a multiplicative
character $\psi$ of $\mathbb{Z}_q^\times$, then
$|\wh{g}|^2=q\delta_{\mathbb{Z}_q\setminus 0}$, i.e., the polynomial 
$$F^{(q)}(G;\mathbf{x})=\prod_{(u,v)\in\overrightarrow{E}}\sum_{b\in\mathbb{Z}_q}\psi(b)x_u^{t-b}x_v^b$$
has $\ell_2$-norm $q^{|E|-|V|}P(G;q)$. (The case $q=3$ is Petersen's graph polynomial reduced modulo $(x_v^3-1:v\in V)$.)
\end{example}

{\small 
\bibliography{fourier} 
\bibliographystyle{plain-annote} 
}

\end{document}